\documentclass[12pt]{amsart}
\usepackage{amsmath,amssymb, epsf}

\newtheorem{theorem}{Theorem}[section]

\newtheorem{definition}[theorem]{Definition}

\newtheorem{lemma}[theorem]{Lemma}
\newtheorem{conj}[theorem]{Conjecture}

\newtheorem{example}[theorem]{Example}
\newenvironment{remark}{\noindent\textbf{Remark}}{}
\newcommand{\three}{\ensuremath{3^{(i)}}}
\newcommand{\thetabh}{\ensuremath{\theta^{BH}}}
\newcommand{\wbar}{\ensuremath{I(\tilde{\alpha})}}
\newcommand{\W}{\ensuremath{{I(\alpha)}}}

\begin{document}

\title[The Hopf algebra of peak functions]{Shifted Quasi-Symmetric Functions
and  the Hopf algebra of peak functions}

\author[N.~Bergeron, S.~Mykytiuk, F.~Sottile, \and S.~van
Willigenburg]{Nantel Bergeron \and Stefan Mykytiuk \and Frank Sottile \and
Stephanie~van~Willigenburg}

\address[Nantel Bergeron and Stefan Mykytiuk and Stephanie van Willigenburg]
        {Department of Mathematics and Statistics\\
        York University\\
        Toronto, Ontario M3J 1P3\\
        CANADA}
\address[Frank Sottile]{Department of Mathematics\\
        University of Wisconsin\\
        Van Vleck Hall\\
        480 Lincoln Drive\\
        Madison, Wisconsin 53706-1388\\
        USA}
\email[Nantel Bergeron]{bergeron@mathstat.yorku.ca}
\urladdr[Nantel Bergeron]{http://www.math.yorku.ca/bergeron}
\email[Stefan Mykytiuk]{Mykytiuk@pascal.math.yorku.ca}
\email[Frank Sottile]{sottile@math.wisc.edu}
\urladdr[Frank Sottile]{http://www.math.wisc.edu/\~{}sottile}
\email[Stephanie van Willigenburg]{steph@pascal.math.yorku.ca}
\urladdr[Stephanie van Willigenburg]%
{http://www-theory.dcs.st-and.ac.uk/\~{}stephvw}

\date{\today}
\thanks{Bergeron supported in part by NSERC}
\thanks{Sottile supported in part by NSF grant DMS-9701755
and NSERC grant  OGP0170279}
\thanks{van Willigenburg is supported in part by the
Leverhulme   Trust}

\maketitle

\def\baselinestretch{.90}
\begin{center}
\begin{minipage}[t]{5in}\footnotesize
{\sc  Summary.}
In his work on $P$-partitions,  Stembridge defined the algebra of peak
functions $\Pi$, which is both a
subalgebra and a retraction of the algebra of quasi-symmetric functions. We
show that $\Pi$ is closed under
coproduct, and therefore a Hopf algebra, and describe the kernel of the
retraction.  Billey and Haiman, in
their work on Schubert polynomials, also defined a new class of
quasi-symmetric functions --- shifted quasi-symmetric functions --- and we
show that
$\Pi$ is strictly contained in the
linear span $\Xi$ of shifted quasi-symmetric functions. We show that $\Xi$
is a coalgebra, and compute the rank of the $n$th graded component.
\bigskip

{\sc R\'esum\'e.}
Dans ses travaux sur les P-partitions, Stembridge d\'efinit l'alg\`ebre
$\Pi$ des fonctions de pics. Cette alg\`ebre peut \^etre vue comme une
sous-alg\`ebre ou un quotient de l'alg\`ebre des fonctions quasi-sym\'etriques.
Nous montrons ici que $\Pi$ est ferm\'ee sous le coproduit,
et est donc une alg\`ebre de Hopf.
Nous d\'ecrivons aussi le noyau du quotient ci-dessus.
D'autre part, dans leurs travaux sur les polyn\^omes de Schubert,
Billey et Haiman ont d\'efini une nouvelle classe de fonctions
quasi-sym\'etriques: les fonctions
quasi-sym\'etrique d\'ecal\'e. Nous montrons que $\Pi$ est strictement
contenue dans l'espace lin\'eaire $\Xi$ des fonctions quasi-sym\'etrique
gauchis. Puis nous montrons que $\Xi$ est une coalg\`ebre et calculons les
dimensions des composantes de degr\'ee $n$.
\end{minipage}
\end{center}

\def\baselinestretch{1.00}

\section{Introduction}

 Schur $Q$ functions first arose in the study of projective representations
of $S_n$ \cite{sch-ur}. Since
then they have appeared  in variety of contexts including the
representations of Lie superalgebras
\cite{serg-eev} and   cohomology classes dual to Schubert cycles in
isotropic Grassmanians \cite{joze-fiak,prag-acz}.
While studying the duality between skew Schur $P$ and $Q$
functions and their
connection to the Schubert calculus of isotropic flag manifolds, we were led
to their
quasi-symmetric analogues: the \emph{peak functions} of  Stembridge  \cite
{stem-bridge}. We show that
\emph{the linear span of peak functions is a Hopf algebra}
(Theorem~\ref{peak-coalg}). We
  also show   that these peak functions are contained in the strictly
larger set of \emph{shifted
quasi-symmetric functions}
 (Theorem~\ref{kw-thetabh}) introduced by Billey and Haiman
\cite{billey-haiman}. We remark that the
quasi-symmetric functions here are not any apparent specialization of the
quasi-symmetric $q$-analogues of Hivert \cite{Hivert}.

{}From extensive calculations, we believe that the set of all shifted
quasi-symmetric functions form a Hopf
algebra, but at present we can only show that:

\emph{The set of all shifted quasi-symmetric functions forms a graded
coalgebra whose $n$th graded component
has rank $\pi_n$, where $\pi_n$ is given by the recurrence}
 \[ \pi_n=\pi_{n-1}+\pi_{n-2} +\pi_{n-4},\]
\emph{with initial conditions $\pi_1=1$, $\pi_2=1$, $\pi_3=2$, $\pi_4=4$.}

We  shall prove this result (Theorems~\ref{sqs-coalg}
and~\ref{graded-rank}) and in addition shall
establish some other properties of these functions.

A composition $\alpha =[\alpha_1,\alpha_2, \ldots ,\alpha_k]$ of a
positive integer $n$ is an ordered list of positive integers whose sum is
$n$. We denote this by $\alpha \vDash n$. We   call the integers
$\alpha_i$ the \textit{components} of $\alpha$, and denote the number of
components in $\alpha$ by $k(\alpha)$. There exists a natural one-to-one
correspondence between compositions of $n$ and subsets of $[n-1]$. If $A=\{
a_1,a_2,\ldots ,a_{k-1}\}\subset[n-1]$, where $a_1<a_2<\ldots <a_{k-1}$,
then $A$ corresponds to the composition, $\alpha =[a_1-a_0,a_2-a_1,\ldots
,a_k-a_{k- 1}]$, where $a_0=0$ and $a_k=n$. For ease of notation, we shall
denote the set corresponding to a given composition $\alpha$ by
$I(\alpha)$. For compositions $\alpha$ and $\beta$ we say that $\alpha$ is
a \emph{refinement} of $\beta$
if $I(\beta)\subset I(\alpha)$, and denote this by
$\alpha\preccurlyeq \beta$.

For any composition $\alpha =[\alpha_1,\alpha_2,\ldots ,\alpha_k]$  we
denote by $M_\alpha$ the
\emph{monomial quasi-symmetric function} \cite{ges-sel}
\[M_\alpha=\sum_{i_1<i_2<\ldots <i_k} x^{\alpha_1}_{i_1}\ldots x^{\alpha
_k}_{i_k}.\]
We define $M_0=1$, where $0$ denotes the unique empty composition of $0$.
We denote by  $F_\alpha$  the \emph{fundamental quasi-symmetric function}
\cite{ges-sel}
\[F_\alpha=\sum_{\alpha\preccurlyeq \beta} M_\beta .\]

\begin{definition}
For any subset $A\subset [n-1]$, let $A+1$ be the subset of $\{2, \dots ,
n\}$ formed from $A$ by
adding $1$ to each element of $A$. Let $\alpha \vDash n$.
Then we define
\[\theta_\alpha =
\sum_{\stackrel{\mbox{\scriptsize $\beta\vDash n$}}{I(\alpha)\subset
I(\beta)\cup I(\beta)+1}}2^{k(\beta)}M_\beta .\]
\label{st-em}\end{definition}
This   is the natural extension of the definition of peak functions given
in~\cite{stem-bridge}.

\begin{example} We shall often omit the brackets that surround the
components of a composition.

If  $\alpha= 21$, then $I(\alpha)=\{ 2\}$, and $I(\alpha)+1=\{3\}$. Hence
$$\theta
_{21}=4M_{21}+4M_{12}+8M_{111}.$$

\end{example}

Let $\Sigma^n$ be the $\mathbb{Z}$-module of quasi-symmetric functions
spanned by
$\{M_\alpha\}_{\alpha\vDash n}$  and let
$\Sigma =\oplus
_{n\geq 0}\Sigma^n$ be the graded $\mathbb{Z}$-algebra of quasi-symmetric
functions.  This is a Hopf
algebra \cite{malv-reut} with coproduct given by
\[\Delta(M_\alpha)=
\sum_{\alpha =\beta \cdot\gamma}M_{\beta }\otimes M_{\gamma},\]
where $\beta\cdot \gamma$ is the concatenation of compositions
$\beta$ and $\gamma$.

\begin{example}
$\Delta (M_{32})=  1\otimes M_{32}+ M_3\otimes M_{2}+
M_{32} \otimes 1$.
\end{example}

We compute the coproduct of the functions $\theta_\alpha$.

\begin{lemma}
For any composition $\alpha\vDash n$ we have that
\begin{eqnarray} \Delta (\theta_\alpha)&=& \sum \theta
_{\epsilon\cdot
a}\otimes\theta_{\phi(b\cdot\zeta)}\label{delta-eq}\end{eqnarray}
where the sum is over all ways of writing $\alpha$ as $\varepsilon \cdot
(a+b) \cdot \zeta$, that is, the
concatenation of compositions $\varepsilon$ and $\zeta$, and a component of
$\alpha$   written as the sum
of numbers $a,b\geq0$. Also $\phi(b\cdot\zeta)= [1+\zeta_1,\zeta_2,\ldots]$
if $b=1$ and $b\cdot\zeta$
otherwise.
\label{mult-coprod}
\end{lemma}

We shall use this result to show that certain subsets of functions
$\theta_\alpha$ span coalgebras (Theorems~\ref{peak-coalg}
and~\ref{sqs-coalg}).

\begin{proof}
Definition~\ref{st-em} is equivalent to
\[
\theta_\alpha =\sum_{\stackrel{\mbox{\scriptsize $\beta\vDash
n$}}{\mbox{\scriptsize $\beta^\ast \preccurlyeq \alpha$}}}
2^{k(\beta)} M_\beta ,
\]
where $\beta^\ast$ is the refinement of $\beta$ obtained by replacing all
components $\beta_i >1$, for
$i>1$, by $[1,\beta_i-1]$. Thus the LHS of equation (\ref{delta-eq}) is
equal to
\begin{eqnarray}
\sum_{\stackrel{
\stackrel{\mbox{\scriptsize$\beta\vDash n$}}
{\mbox{\scriptsize$\beta^\ast\preccurlyeq\alpha$}}}
{\mbox{\scriptsize$\beta=\gamma\cdot\delta$}}} 2^{k(\beta)}
M_\gamma \otimes M_\delta
&=&\sum_{\stackrel{\mbox{\scriptsize$\gamma\cdot\delta\vDash
n$}}{(\gamma\cdot\delta)^\ast\preccurlyeq\alpha}}
 2^{k(\gamma)}M_\gamma\otimes 2^{k(\delta)}M_\delta.
\label{2delta-eq}
\end{eqnarray}

Let $2^{k(\gamma)}M_\gamma\otimes 2^{k(\delta)}M_\delta$ be a term of this
sum, with
$\gamma\vDash m$. This term can only appear in one summand on the RHS of
equation (\ref{delta-eq}),
namely $\theta_{\varepsilon\cdot a}\otimes \theta_{\phi(b\cdot\zeta)}$ with
$\varepsilon\cdot a\vDash
m$.
To show that it does indeed appear, we need to prove that $\gamma^\ast
\preccurlyeq
\varepsilon\cdot a$ and $\delta^\ast \preccurlyeq \phi(b\cdot\zeta)$. Let
$\delta^{\ast\ast}$ be the
refinement of $\delta^\ast$ obtained by replacing the part $\delta_1$ by
$[1,\delta_1-1]$ if $\delta_1>1$.
We have that
\[\gamma^\ast\cdot \delta^{\ast\ast}=(\gamma\cdot\delta)^\ast\preccurlyeq
\varepsilon\cdot(a+b)\cdot\zeta
,\]
which implies that $\gamma^\ast \preccurlyeq \varepsilon\cdot a$, and
$\delta^{\ast\ast}\preccurlyeq
b\cdot\zeta\preccurlyeq \phi(b\cdot \zeta)$.

If $\delta_1=1$ then $\delta^\ast =\delta^{\ast\ast}\preccurlyeq\phi
(b\cdot \zeta)$. However, if $\delta
_1>1$ then there are two possible cases: either $\delta_1\leq b$, or $b=1$
and $\delta_1-1\leq\zeta_1$. In
the former case $\delta^\ast \preccurlyeq b\cdot\zeta =\phi (b\cdot
\zeta)$, while in the latter, $\delta
_1\preccurlyeq 1+\zeta_1$, whence $\delta^\ast \preccurlyeq
[1+\zeta_1,\zeta_2,\ldots ]=
\phi ( b\cdot\zeta)$.

Conversely, let $ 2^{k(\gamma)}M_\gamma\otimes 2^{k(\delta)} M_\delta$ be a
term belonging to a
tensor $\theta_{\varepsilon\cdot a}\otimes \theta_{\phi(b\cdot \zeta)}$ on
the RHS of equation
(\ref{delta-eq}). To show that it appears in equation  (\ref{2delta-eq}) we
must prove that $(\gamma\cdot
\delta )^\ast\preccurlyeq \varepsilon\cdot (a+b)\cdot \zeta$. We have that
$\gamma^\ast\preccurlyeq
\varepsilon\cdot a$ and $\delta^\ast\preccurlyeq \phi(b\cdot \zeta)$, which
imply that
\[(\gamma\cdot \delta)^\ast= \gamma^\ast\cdot\delta^{\ast\ast}\preccurlyeq
\gamma^\ast\cdot\delta^\ast
\preccurlyeq \varepsilon\cdot a\cdot \phi(b\cdot \zeta).\]

If $b>1$ then
$$(\gamma\cdot\delta)^\ast\preccurlyeq \varepsilon\cdot a\cdot \phi(b\cdot
\zeta)=\varepsilon\cdot a\cdot
b\cdot \zeta\preccurlyeq \varepsilon\cdot (a+b)\cdot\zeta.$$
If $b=1$ then $\delta^\ast \preccurlyeq
\phi(b\cdot\zeta)=[1+\zeta_1,\zeta_2,\ldots ]$ implies that
\[\delta^{\ast\ast}=[1, \ldots ]\preccurlyeq [1,\zeta_1,\ldots ]=b\cdot
\zeta .\]
Therefore,
\[(\gamma\cdot\delta)^\ast=\gamma^\ast\cdot\delta^{\ast\ast}\preccurlyeq
\varepsilon\cdot a\cdot
b\cdot\zeta\preccurlyeq \varepsilon\cdot (a+b)\cdot \zeta\]
as desired.
\end{proof}

\section{The peak Hopf algebra}

\begin{definition}
For any composition $\alpha =[\alpha_1,\alpha_2,\ldots ,\alpha_k]$ we say
that $\theta_\alpha$ is a
\emph{peak function} if $\alpha_i = 1\Rightarrow i=k$.\end{definition}

Observe that if $\theta_\alpha$ is a peak function and $\alpha\vDash n$,
then $I(\alpha)\subset \{2,\ldots ,n-1\}$ such that no two
$i$ in $I(\alpha)$ are consecutive.

Let  $\Pi^n$ be the $\mathbb{Z}$-module spanned by
all peak functions $\theta
_\alpha$, $\alpha\vDash n$, and let $\Pi=\oplus_{n\geq 0}\Pi^n$. This was
studied by Stembridge
\cite{stem-bridge} who
showed that the peak functions are F-positive, are closed under product,
and form a basis for $\Pi$, and so
the rank of   $\Pi^n$ is  the
$n$th Fibonacci number. In addition we also
know the following about the \emph{algebra of peaks}, $\Pi$.

\begin{theorem} $\Pi$ is closed under coproduct.
\label{peak-coalg}\end{theorem}

\begin{proof}
If all components of a composition $\alpha$, except perhaps the last, are
greater than $1$, then the same is
true for all compositions $\varepsilon\cdot a$ and $\phi(b\cdot\zeta)$
appearing in the RHS of equation
(\ref{delta-eq}).
\end{proof}

Let $\Theta$ be the $\mathbb{Z}$-linear map from $\Sigma $ to $\Pi $
defined by
$\Theta(F_\alpha)=\theta_{\Lambda(\alpha)}$,
where $\Lambda(\alpha)$ is
the composition formed from $\alpha =[\alpha_1,\alpha_2,\ldots ,\alpha
_k]$ by adding together
adjacent components $\alpha_i ,\alpha_{i+1},\ldots , \alpha_{i+j}$ where
$\alpha_{i+l}=1$ for
$l=0,\ldots , j-1$, and either $\alpha_{i+j}\neq 1$, or
  $i+j=k$.

\begin{example}
If $\alpha=31125111$ then $\Lambda (\alpha)=3453$.
\end{example}

Stembridge \cite{stem-bridge} showed that $\Theta:\Sigma \rightarrow \Pi$ is
a graded surjective ring
homomorphism, and was an analogue of the retraction from the algebra of
symmetric functions to
Schur $Q$ functions. It is clear from our proof above that this morphism is
in fact a Hopf homomorphism.
We can
describe the kernel of $\Theta$ as follows.

\begin{lemma}
The non-zero differences $F_\alpha-F_{\Lambda (\alpha)}$ form a basis of
the kernel of $\Theta$.
\end{lemma}

\begin{proof}
Each difference $F_\alpha -F_{\Lambda(\alpha)}$ is in the kernel of
$\Theta$ as $\Theta(F_\alpha -F_{\Lambda(\alpha)})=0$ since
$\Lambda(\Lambda(\alpha))=\Lambda(\alpha)$.
In addition, the
non-zero differences are linearly independent as they have different leading
terms. Letting $f_n$ denote the
$n$th Fibonacci number, there are  $2^{n-1}-f_n$
such differences, and since
 \begin{eqnarray*}
\dim\ker\Theta &=& \dim \Sigma^n - \dim \Pi^n\\
&=&2^{n-1} -f_n,
\end{eqnarray*}
our result follows.
\end{proof}

\section{The   coalgebra of shifted quasi-symmetric functions}
\begin{definition}
For any composition $\alpha =[\alpha_1,\alpha_2,\ldots ,\alpha_k]\vDash n$
we say that
$\theta_\alpha$ is a
\emph{shifted quasi-symmetric function} (sqs-function) if $n\le 1$ or
$\alpha_1> 1$.\end{definition}

Observe that if $\theta_\alpha$ is an sqs-function  and $\alpha\vDash n$,
then $I(\alpha)\subset
\{2,\ldots  ,n-1\}$.

For integers $n\geq 0$, let $\Xi^n$ be the $\mathbb{Z}$-module spanned by
all sqs-functions
$\theta_\alpha$, $\alpha\vDash n$, and let $\Xi=\oplus_{n\geq 0}\Xi^n$.

\begin{theorem}
$\Xi$ is closed under coproduct.
\label{sqs-coalg}
\end{theorem}

\begin{proof}
If  the first component of a composition $\alpha$ is greater than $1$, then
the same is true for all
compositions $\varepsilon\cdot a$ and $\phi(b\cdot\zeta)$ appearing in the
RHS of equation
(\ref{delta-eq}).
\end{proof}

Unlike   peak
functions \cite{stem-bridge},   sqs-functions are not
$F$-positive since
$$\theta_{211}=F_{22}+ F_{112}+
2F_{121}+
F_{211}-
F_{1111}.$$
 \begin{definition}
For any composition, $\alpha\vDash n$, we define the complement $\alpha^c$
of $\alpha$ to be the
composition for which $I(\alpha^c) = (I(\alpha))^c$, the set complement of
$I(\alpha)$ in $[n-1]$. We
define the graph
$G(\alpha)$ of $\alpha$ to be the graph obtained from

$$
\epsfxsize=1.8in\epsfbox{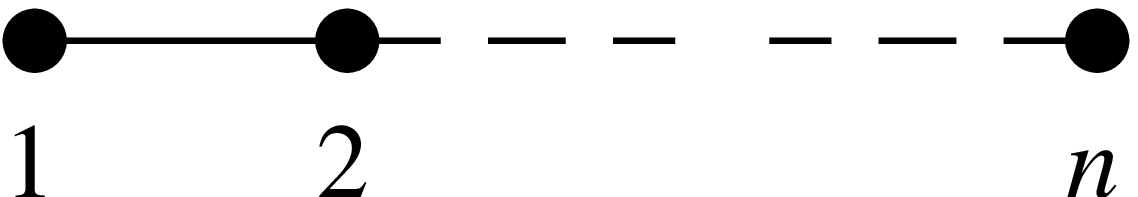}
$$

by removing the edge $(i,i+1)$ if and only if $i\in I(\alpha)$.
\end{definition}

Observe that $G(\alpha^c)$ contains the edge $(i,i+1)$ if and only if this
edge is not contained in
$G(\alpha)$. These graphs will be used later to simplify the proof of
Theorem~\ref{kw-thetabh}.

Let a \emph{word} of length $n$ be any $n$-tuple, $w_1w_2\ldots w_n$, and
let a
\emph{binary word} of length $n$ be a word $w_1w_2\ldots w_n$ wuch that
$w_i\in \{ 0,1\}$ for all
$i$.   For $2\leq
i\leq n-1$, let us denote by  $\three $ the composition $
[1^ {i-2}, 3,1^{n-i-1}]$ of
$n $.
For some subset $S\subset \{2, \ldots ,n-1\}$, let us denote by
$\bigwedge_{i\in S}\three$ the
composition of $n$ for which $G(\bigwedge_{i\in S}\three)$ has an edge
between vertices $i$ and $i+1$ if
and only if an edge exists between vertices $i$ and $i+1$ in $G(\three)$ for
some $i\in S$.

\begin{example}
Let $S=\{ 2,3\}\subset [3]$. Then $G(3^{(2)})$ is
$$
\epsfxsize=1.55in\epsfbox{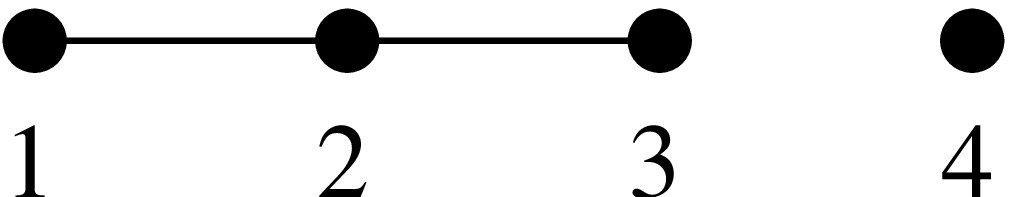}
$$
and  $G(3^{(3)})$ is
$$
\epsfxsize=1.55in\epsfbox{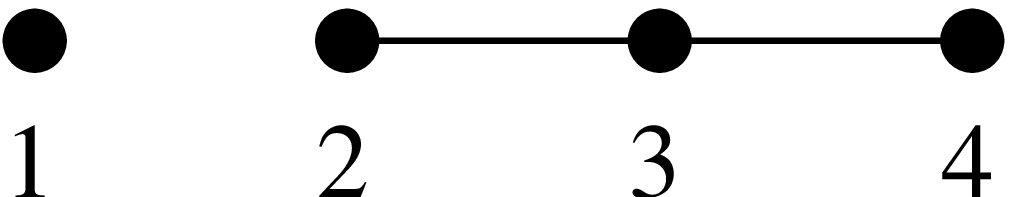}
$$
hence $G(\bigwedge_{i\in S}\three)$ is
$$
\epsfxsize=1.55in\epsfbox{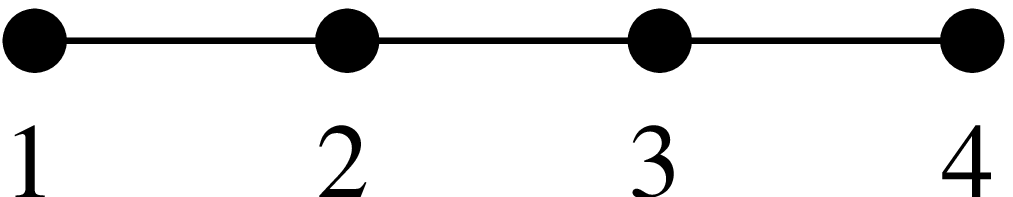}
$$
so $\bigwedge_{i\in S}\three$ is the composition $4$.
\end{example}

\begin{definition}\cite{billey-haiman}
Let $\alpha$ be a composition of $n$. Let ${\mathcal A}(I(\alpha))$ denote
the set of
all sequences $j_1\leq
j_2\leq\ldots \leq j_n$ in $\mathbb{N}$ such that we do not have
$j_{i-1}=j_i=j_{i+1}$ for any $i\in
I(\alpha)$. The \emph{shifted quasi-symmetric function} $\thetabh_\alpha$  is
given by
\[\thetabh_\alpha =
\sum_{\stackrel{\mbox{\scriptsize$J=(j_1,\ldots, j_n)$}}{
\stackrel{\mbox{\scriptsize$j_1\leq\ldots\leq j_n$}}{
J\in{\mathcal A}(I(\alpha))}}} 2^{| j|}x_{j_1}\ldots
x_{j_n},\]where $|j|$ denotes the number of distinct values $j_i$ in
$J$.\label{bi-ha}
\end{definition}

\begin{theorem}
For any sqs-function $\theta_\alpha$ we have that  $\theta_\alpha=\thetabh
_\alpha$.\label{kw-thetabh}
\end{theorem}

\begin{proof}
For each   $i\in \W\subset[n-1]$, $j_{i-1}=j_i=j_{i+1}$ is forbidden in
any monomial
\[ x_{j_1}x_{j_2}\ldots x_{j_i}\ldots x_{j_n}\]
appearing as a summand of the function $\thetabh_\alpha$. This is
equivalent to saying that $M_\beta$ is
a
summand of $\thetabh_\alpha$ if and only if $G(\three)\not\subset
G(\beta)$ for all $i\in \W$.
Therefore at least one
of $i-1$ or $i$ must be   the largest label of a vertex in a connected
component in $G(\beta)$.

Now when going from compositions of $n$ to subsets of $[n-1]$ we can do so
using our graphs, $G$. All
we have to do is list the label of the vertex that is the largest in each
connect component, not listing $n$. We
call these vertices the \emph{end-points}.
We are now in a position to prove the equivalence of    Definitions
~\ref{st-em} and
~\ref{bi-ha} for sqs-functions.

The powers of 2 agree  so we need only show that the indices of summation
do too. To see this, take any
sqs-function $\theta_\alpha$ and let $i\in\W$. Then
$M_\beta$ is a summand in $\thetabh_\alpha$    if at least one of $i-1$ or
$i$ is an end-point in $G(\beta)$.
Therefore $i$ or $i-1$ belongs to $I(\beta)$, and
$M_\beta$  is a summand of $\theta_\alpha$. Conversely, if $M_\beta $  is
a summand of $\theta
_\alpha$,   then this implies that for
each $i\in \W$, we have that $i-1$ or $i$ belongs to $I(\beta)$, so one of
$i-1$ or $i$
is an end-point in $G(\beta)$,
so $M_\beta$ is a summand of $\thetabh_\alpha$.
\end{proof}

 \section{A basis for $\Xi$}

\begin{definition}
Let $\theta_\alpha$ be an sqs-function and $\alpha \vDash n$. We define an
internal peak
$i\in \W$ such that $i-
1,i+1\not\in \W$, and $i\in\{ 3,\ldots
n-2\}$.
\end{definition}

\begin{remark}
Observe that the occurrence of an internal peak in the $i$th position in
$I(\alpha)=\{ w_1,w_2,\ldots\}$,  where
$w_1<w_2<\ldots$, is equivalent to having two components of $\alpha$, say
$\alpha_i,\alpha_{i+1}$
such that $\alpha_{i+1}\geq 2$, and $\alpha_i\geq 2$ if $i\neq1$, or
$\alpha_i\geq 3 $ if $i=1$.
\end{remark}

We can now describe the basis of $\Xi$ as follows.

\begin{theorem}
The coalgebra  $\Xi$ has a basis consisting of all sqs-functions
$\theta_\alpha$   where  $\W$ contains no
internal peak.\label{xi-basis}
\end{theorem}

We sketch the proof of Theorem 4.2 later.

\begin{theorem}
The rank of $\Xi^n$ is given by the recurrence
\[ \pi_n=\pi_{n-1}+\pi_{n-2} +\pi_{n-4},\]
with initial conditions $\pi_1=1$, $\pi_2=1$, $\pi_3=2$,
$\pi_4=4$.\label{graded-rank}
\end{theorem}

This recurrence was suggested by a superseeker query \cite{sloan-seeker}.

\begin{proof}
By direct calculation we obtain that  $\pi_1=1$, $\pi_2=1$, $\pi_3=2$,
and $\pi_4=4$.

To obtain our recurrence, we observe that for each sqs-function,
$\theta_\alpha$ where $\alpha\vDash n$,
we can encode $I(\alpha)$
as  a binary word of length
$n-2$, by
placing a $1$ in position $i-1$ if $i$ is contained in $I(\alpha)$, and $0$
otherwise. By this
one-to-one
correspondence we see that $I(\alpha)$  contains no internal peak
if its corresponding binary word
does not
contain $010$ as a subword.

We therefore count binary words of length $n$ that avoid the subword
$010$. Appending either $1$ or $0$
to such a binary word of length $n-1$ gives one of length $n$, provided that
we have not created the
subword $010$ in the last three positions. Let $a_n$, $b_n$, $c_n$, and
$d_n$ enumerate those binary
words of length  $n-2$ that avoid the subword $010$ and end in,
respectively $00$, $01$, $10$, and $11$.
We then obtain the following 4 simultaneous recursions.
\[a_n=a_{n-1} +c_{n-1}, \ b_n=a_{n-1} +c_{n-1}, \  c_n=d_{n-1}, \ d_n=b_{n-1}
+d_{n-1}. \]

Clearly the number of $I(\alpha)$s in $[n-1]$ with no internal peaks is
given by
\[ \pi_n=a_n+b_n+c_n+d_n,\]
However by substituting in our recurrences we obtain
\begin{eqnarray*}
\pi_n&=& a_n+b_n+c_n+d_n\\
&=& 2a_{n-1}+b_{n-1}+2c_{n-1}+2d_{n-1}\\
&=& \pi_{n-1}+a_{n-1}+c_{n-1}+d_{n-1}\\
&=& \pi_{n-1}+a_{n-2}+b_{n-2}+c_{n-2}+2d_{n-2}\\
&=& \pi_{n-1}+\pi_{n-2} +d_{n-2}\\
&=& \pi_{n-1}+\pi_{n-2} +b_{n-3} +d_{n-3}\\
&=& \pi_{n-1}+\pi_{n-2} +a_{n-4}+b_{n-4} +c_{n-4}+d_{n-4}\\
&=& \pi_{n-1}+\pi_{n-2} +\pi_{n-4}. \\
\end{eqnarray*}
\end{proof}

We say that $M_\beta$ is a maximal term of $\theta_\alpha$ if for any
$\gamma$ higher in the partial order
of compositions $M_\gamma$ is not a summand of $\theta_\alpha$. The
following lemma is stated without
proof.

\begin{lemma}
Let $\theta_\alpha$ be an sqs-function. Consider the collection $S$ of all
possible sets derived from
$I(\alpha)$   by
adding either $i-1$ or
$i+1$ to
$\W$ for all internal peaks $i\in \W$. If $M_\beta$ is a maximal term of
$\theta_\alpha$, then $\beta$
is
derived from
\[\bigwedge_{\stackrel{\mbox{\scriptsize$i\in (\wbar )^c$}}{\wbar\in S}}
\three \]
by adding adjacent components equal to 1 together to give a component equal
to 2 as often as possible.
\label{big-term}
\end{lemma}

 \begin{lemma} Let $\theta_\alpha$ be an sqs-function, and let $ \W$
have an internal peak in the $j$th   position, then we have the following
linear relation
\begin{eqnarray*}\theta_\alpha &=&\theta_{[\alpha_1,\ldots
,\alpha_j-1,1,\alpha_{j+1},\ldots
,\alpha_k]}+ \theta_{[\alpha_1,\ldots ,\alpha_j,1,\alpha_{j+1}-1,\ldots
,\alpha_k]}\\
&&- \theta_{[\alpha_1,\ldots ,\alpha_j-1,1,1,\alpha_{j+1}-1,\ldots ,\alpha
_k]}.\end{eqnarray*}\label{rel-ation}
\end{lemma}

\begin{proof}
By Definition~\ref{bi-ha} we have that the leading terms of $\theta
_\alpha$ determine the other summands that belong to $\theta_\alpha$.
Hence by Lemma~\ref{big-term} it follows that the summands of $\theta
_\alpha$ will be the union of the summands of $ \theta_{[\alpha_1,\ldots
,\alpha_j-1,1,\alpha_{j+1},\ldots ,\alpha_k]}$ and $\theta_{[\alpha
_1,\ldots ,\alpha_j,1,\alpha_{j+1}-1,\ldots ,\alpha_k]}$.  However, those
summands that appear in both will be duplicated. By definition these will be
the summands of $\theta_{[\alpha_1,\ldots ,\alpha_j-1,1,1,\alpha
_{j+1}-1,\ldots ,\alpha_k]}$, and the result follows.
\end{proof}

\noindent\textit{Sketch of proof of Theorem~\ref {xi-basis}.} {}From our
relation in Lemma~\ref{rel-ation}, it follows that any $\theta_\alpha$ can
be rewritten as a
linear combination of functions $\theta_{\tilde{\alpha}}$, where
$I(\tilde{\alpha})$ contains no internal
peaks. In addition, by Lemma~\ref{big-term} and definition~\ref{bi-ha} we
have that the set of all
sqs-functions $\theta_\alpha$ where $I(\alpha)$ contains no internal peaks
is linearly independent and thus
form a basis for $\Xi$.\qed


\begin{thebibliography}{99}
\bibitem{billey-haiman}  S. Billey, and M. Haiman, \textit{ Schubert
  polynomials for the classical groups},
  Journal of AMS, \textbf{8} (1995) 443--482.
\bibitem{ges-sel}
  I. Gessel, \textit{Multipartite $P$-partitions and inner products of skew
  {S}chur  functions},  Contemporary Mathematics, \textbf{34} (1984) 289--301.

\bibitem{Hivert}
F. Hivert, \textit{Quasi-symmetric functions and Hecke algebra action},
Proceeding of the 10th
FPSAC, Toronto (1998) 377--388.

\bibitem{joze-fiak} T. J\'{o}zefiak, \textit{ Schur $Q$ functions and
  cohomology of isotropic Grassmanians},
  Math.~Proc. Camb.~Phil.~Soc., \textbf{109} (1991) 471--478.
\bibitem{malv-reut}  C. Malvenuto and C. Reutenauer, \textit{ Duality
  between quasi-symmetric functions
  and the {S}olomon descent algebra},
  Journal of  Algebra, \textbf{177} (1995) 967--982.
\bibitem{sloan-seeker}
     N.J.A.~Sloane,\textit{An on-line version of the
      Encyclopedia of integer sequences},
      Elect.~J.~Combin., \textbf{1} No.1,  (1994)  F1.

      {\tt http://akpublic.research.att.com/\~{}njas/sequences/ol.html}.
\bibitem{prag-acz}  P. Pragacz, \textit{Algebro-geometric applications of
  {S}chur  $S$- and $Q$- polynomials},
  Topics in Invariant Thoery, S\'{e}minaire d'Alg\`{e}bre Dubreil-Malliavin
   1989-90, Springer-Verlag,
   (1991) 130--191.
\bibitem{sch-ur}  I. Schur, \textit{ \"{U}ber die Darstellung der
 symmetrischen und der alternierenden
 Gruppe durch gebrochene lineare Substitutionen},
 J.~reine angew. Math. , \textbf{139} (1911) 155--250.
\bibitem{serg-eev}  A.N. Sergeev, \textit{ The tensor algebra of the
 identity representation as a module
 over the Lie superalgebras $gl(n,m)$ and $Q(n)$},
 Math. USSR Sbornik, \textbf{51} (1985) 419--427.
\bibitem{stem-bridge}  J. Stembridge, \textit{ Enriched{ P}-partitions},
 Trans. AMS, \textbf{349} (1997) 763--788.
\end{thebibliography}
\end{document}